\documentclass{article}

\usepackage{amssymb}
\usepackage{appendix}
\usepackage{latexsym}
\usepackage{mathrsfs}
\usepackage{extarrows}
\usepackage[dvips]{graphicx}

\newtheorem{theo}{Theorem}[section]

\newtheorem{remark}[theo]{Remark}
\newtheorem{ques}{Question}
\newtheorem{lemma}[theo]{Lemma}
\newtheorem{claim}[theo]{Claim}

\newtheorem{con}[theo]{Conjecture}
\newtheorem{prop}[theo]{Proposition}
\newtheorem{fact}[theo]{Fact}
\newtheorem{defi}[theo]{Definition}

\def\q{\hspace*{\fill}$\Box$\medskip}

\def\endproofbox{\hskip 1.3em\hfill\rule{6pt}{6pt}}

{%
\noindent{\it Proof}
}%
{%
 \quad\hfill\endproofbox\vspace*{2ex}
}

\begin{document}
%\date{}
\title{Lagrangian densities of  some $3$-uniform hypergraphs}
\author{Zilong Yan \thanks{School of Mathematics, Hunan University, Changsha 410082, P.R. China. Email: zilongyan@hnu.edu.cn.} \and Yuejian Peng \thanks{ Corresponding author. School of Mathematics, Hunan University, Changsha, 410082, P.R. China. Email: ypeng1@hnu.edu.cn. \ Supported in part by National Natural Science Foundation of China (No. 11931002 and No. 11671124).}
}

\maketitle

\begin{abstract}
The Lagrangian density of an $r$-uniform  hypergraph $H$ is $r!$ multiplying the supremum of the Lagrangians of all $H$-free $r$-uniform  hypergraphs.
For an $r$-uniform graph $H$ with $t$ vertices, it is clear that $\pi_{\lambda}(H)\ge r!\lambda{(K_{t-1}^r)}$. We say that an $r$-uniform hypergraph $H$ with $t$ vertices is $\lambda$-perfect if $\pi_{\lambda}(H)= r!\lambda{(K_{t-1}^r)}$.
A theorem of Motzkin and Straus implies that all  $2$-uniform graphs are $\lambda$-perfect.  It is interesting to understand what kind of hypergraphs are $\lambda$-perfect. The property `$\lambda$-perfect' is monotone in the sense that an $r$-graph obtained by removing an edge from a  $\lambda$-perfect $r$-graph (keep the same vertex set) is $\lambda$-perfect. It's interesting to understand the relation between  the number of edges in a hypergraph and the `$\lambda$-perfect' property. We propose that the number of edges in a hypergraph no more than the number of edges in a linear hyperpath would guarantee the `$\lambda$-perfect' property. We show some partial result to support this conjecture. We also give some partial result to support the conjecture that the disjoint union of two $\lambda$-perfect $r$-uniform hypergraph is $\lambda$-perfect. We show that the disjoint union of a $\lambda$-perfect $3$-graph and $S_{2,t}=\{123,124,125,126,...,12(t+2)\}$ is perfect. This result implies the earlier result of Heftz and Keevash, Jiang, Peng and Wu, and several other earlier results.

Let $\Pi_{r}=\{\pi({\mathcal F}):{\mathcal F} {\rm \ is \  a \ family \ of\ } r-{\rm uniform \ graphs}. \}$, where $\pi({\mathcal F})$ is the Tur\'{a}n density of ${\mathcal F}$.  Erd\H{o}s-Stone-Simonovits showed  that $\Pi_{2}=\{0, {1 \over 2}, {2 \over 3}, {3 \over 4}, \ldots \}$.
However, very few  are known for $\Pi_{r}$ when $r\ge 3$.  Sidorenko showed that the Lagrangian density of  an $r$-uniform graph  is the same as  the  Tur\'{a}n density of its extension. In particular, these two densities of $H$ equal if every pair of vertices of  $H$  is contained in an edge. Sidorenko  gave  the Lagrangian densities of  infinitely many hypergraphs and this led to first know infinitely many numbers in $\Pi_{r}$ for $r\ge 3$. Determining the Tur\'an density of $K_4^{3-}$ (a well-known conjecture of Frankl and F\"uredi) is one of the major open problems in hypergraph Tur\'an problems. The result of Sidorenko implies that the Lagrangian density and the Tur\'an density of $K_4^{3-}$  equal. It would be very interesting to obtain the  Tur\'an density of $K_4^{3-}$  by determining the Lagrangian density of $K_4^{3-}$. As a first attempt, we determine the Lagrangian density of $K_4^{3-}\sqcup e$, where $e$ is a single edge. We also find an irrational Lagrangian density. The advantage of transferring to the Lagrangian density is that we can assume that $G$ is `dense' (so covering pairs) when considering the maximum Lagrangian of an $H$-free $r$-graph $G$. This assumption makes the structural analysis `nicer' in some cases, so we hope that this method helps us  to better understand the set $\Pi_r$.

\end{abstract}


\section{Introduction}
\subsection{ Notations and definitions}
For a set $V$ and a positive integer $r$, let $V^ r$ denote the family of all $r$-subsets of $V$. An {\em $r$-uniform graph} or {\em $r$-graph $G$} consists of a set $V(G)$ of vertices and a set $E(G) \subseteq V(G) ^r$ of edges. Let $|G|$ denote the number of edges of $G$. An edge $e=\{a_1, a_2, \ldots, a_r\}$ will be simply denoted by $a_1a_2 \ldots a_r$. An $r$-graph $H$ is a {\it subgraph} of an $r$-graph $G$, denoted by $H\subseteq G$, if $V(H)\subseteq V(G)$ and $E(H)\subseteq E(G)$. A subgraph of $G$ {\em induced} by $V'\subseteq V$, denoted as $G[V']$, is the $r$-graph with vertex set $V'$ and edge set $E'=\{e\in E(G):e \subseteq V'\}$. For $S\subseteq V(G)$, let $G-S$ denote the subgraph of $G$ induced by $V(G)\setminus S$. Let $K^{r}_t$ denote the complete $r$-graph on $t$ vertices. Let $K^{r-}_t$ be obtained by removing one edge form the complete $r$-graph on $t$ vertices. For a positive integer $n$,  let $[n]$ denote $\{1, 2, 3, \ldots, n\}$.

Given an $r$-graph $F$,  an $r$-graph $G$ is called $F$-free if it does not contain an isomorphic copy of $F$. For a fixed positive integer $n$ and an $r$-graph $F$, the {\em Tur\'an number} of $F$, denoted by $ex(n,F)$, is the maximum number of edges in an $F$-free $r$-graph on $n$ vertices.
 An averaging argument of Katona, Nemetz and Simonovits \cite{KNS} shows that the sequence ${ ex(n,F) \over {n \choose r } }$ is  non-increasing. Hence $\lim_{n\rightarrow\infty} { ex(n,F) \over {n \choose r } }$ exists. The {\em Tur\'{a}n density} of $F$ is defined as $$\pi(F)=\lim_{n\rightarrow\infty} { ex(n,F) \over {n \choose r } }.$$ Denote $$\Pi_{r}=\{\pi({\mathcal F}):{\mathcal F} {\rm \ is \  a \ family \ of\ } r-{\rm uniform \ graphs}. \}.$$

For 2-graphs, Erd\H{o}s-Stone-Simonovits determined the Tur\'an numbers of all non-bipartite graphs asymptotically. Their Result implies that
$\Pi_{2}=\{0, {1 \over 2}, {2 \over 3}, {3 \over 4}, \ldots \}$.
Very few results are known for $r\ge 3$ and a  survey on this topic can be found in Keevash's survey paper \cite{Keevash}.  The hypergraph Lagrangian method has been helpful in  hypergraph extremal problems.

\begin{defi}
Let $G$ be an $r$-graph on $[n]$ and let
  $\vec{x}=(x_1,\ldots,x_n) \in [0,\infty)^n$. Define the {\em Lagrangian} function
$$\lambda (G,\vec{x})=\sum_{e \in E(G)}\prod\limits_{i\in e}x_{i}.$$
\end{defi}
The {\em Lagrangian} of
$G$, denoted by $\lambda (G)$, is defined as
 $$\lambda (G) = \max \{\lambda (G, \vec{x}): \vec{x} \in \Delta \},$$
where $$\Delta=\{\vec{x}=(x_1,x_2,\ldots ,x_n) \in [0, 1])^{n}: x_1+x_2+\dots+x_n =1 \}.$$

The value $x_i$ is called the {\em weight} of the vertex $i$ and a vector $\vec{x} \in {\Delta}$ is called a {\em feasible weight vector} on $G$.
A feasible weight vector  $\vec{y}\in {\Delta}$ is called an {\em optimum weight vector} for $G$ if $\lambda (G, \vec{y})=\lambda(G)$.

In \cite{MS}, Motzkin and Straus established a connection between the Lagrangian of a $2$-graph and it's maximum complete subgraphs.
\begin{theo} {\em(\cite{MS})} \label{MStheo}
If $G$ is a $2$-graph in which a maximum complete subgraph has  $t$ vertices, then
$\lambda(G)=\lambda(K_t^2)={1 \over 2}(1 - {1 \over t})$.
\end{theo}

 They also applied this connection to give another proof of the theorem of Tur\'an  on the Tur\'an density of complete graphs. Since then the Lagrangian method has been a useful tool in hypergraph extremal problems. Earlier applications include that Frankl and R\"{o}dl \cite{FR} applied it in disproving the long standing jumping constant conjecture of Erd\H{o}s. Sidorenko \cite{Sidorenko-89} applied Lagrangians of hypergraphs to first find infinitely many Tur\'an densities of  hypergraphs. More recent developments of this method were obtained in  \cite{Pikhurko, BT,  HK, NY, BIJ, NY2, Jenssen, JPW, HPW, YP}.
 In addition to its applications in extremal problems,  it is interesting in its own right to determine the  maximum Lagrangian of $r$-graphs with certain properties as remarked by Hefetz and Keevash \cite{HK}. For example, an interesting conjecture of Frankl-F\"uredi \cite{FF} states that the the maximum  Lagrangian among all $r$-graphs with $m$ edges achieves on the $r$-graph whose edges are the first $r$-tuples in colex order.  Talbot \cite{T} made a first breakthrough in  confirming  this conjecture for some cases. Subsequent progress  were made in \cite{Tyo, LLP2018, TPZZ2}. Recently, Gruslys-Letzter-Morrison \cite{GLM2018} confirmed this conjecture for $r=3$ and sufficiently large $m$, and showed that the conjecture is not always true for $r\ge 4$.  We focus on  the Lagrangian density of an $r$-graph in this paper.

The {\em Lagrangian density } $\pi_{\lambda}(F)$ of  an $r$-graph $F$ is defined to be
$$\pi_{\lambda}(F)=\sup \{r! \lambda(G): G \:\: {\rm is} \:\: F{\text-}{\rm free}\}.$$
A pair of vertices $\{i, j\}$ is {\em covered}  in a hypergraph $F$ if there exists  an edge $e$ in $F$ such that  $\{i, j\}\subseteq e$. We say that $F$ covers pairs if every pair of vertices in $F$ is covered.
 Let $r\ge 3$ and $F$ be an $r$-graph.  The  {\em extension} of $F$, denoted by $H^F$ is
obtained as follows: For each pair of vertices $v_i$ and $v_j$ not covered in $F$, we add a set $B_{ij}$ of $r-2$ new vertices and the edge $\{v_i,v_j\} \cup B_{ij}$, where the $B_{ij}$'s are pairwise disjoint over all such pairs $\{i,j\}$.

The Lagrangian density  is closely related to the Tur\'an density.
The following proposition is implied by Theorem $2.6$ in \cite{Sidorenko-87} (see Proposition $5.6$ in  \cite{BIJ} and Corollary $1.8$ in \cite{Sidorenko-89} for the explicit statement). %It gives the  connection between the Lagrangian density of a hypergraph  $F$ and the Tur\'an density of its extension.
\begin{prop}{\em (\cite{Sidorenko-87, BIJ,Sidorenko-89})}\label{relationlt} Let $F$ be an $r$-graph. Then \\
$(i)$ $\pi(F)\le \pi_{\lambda}(F);$ \\
$(ii)$ $\pi(H^F)=\pi_{\lambda}(F).$ In particular, if $F$ covers pairs, then $\pi(F)= \pi_{\lambda}(F).$
\end{prop}

 The Lagrangian density of the enlargement of a tree satisfying Erd\H{o}s-Sos's conjecture is determined by Sidorenko \cite{Sidorenko-89} and Brandt-Irwin-Jiang \cite{BIJ}. This  first gave  the Tur\'an densities of  infinitely many hypergraphs.
%{\color{red} As a possible generalization of the problem of computing $ex(n,K_{3})$ to 3-graphs, Keevash and Mubayi \cite{KM} 2004 showed that all cancellative 3-graphs on $n$ with close to the maximum number of edges must be structurally close to the complete balanced 3-partite 3-graph. }
Pikhurko \cite{Pikhurko} determined the Lagrangian density of a $4$-uniform tight  path of length 2  and this led  to confirm the conjecture of Frankl-F\"uredi on the  Tur\'an number of its extension, the $r$-uniform genearlized triangle for the case $r=4$.   Norin and Yepremyan \cite{NY2} determined  for $r=5$ or $6$  by extending the earlier result of Frankl-F\"uredi in \cite{FF}. Jenssen \cite{Jenssen} determined the Lagrangian density of a path of length 2 formed by two edges intersecting at   $r-2$ vertices for $r=3, 4, 5, 6, 7$.
Hefetz and Keevash \cite{HK} determined the Lagrangian density of a $3$-uniform matching  of size 2. Jiang-Peng-Wu \cite{JPW} obtained for any $3$-uniform matching. The case for  an $r$-uniform matching  of size 2   was given  in \cite{NWY} (independently, in \cite{WPC} for $r=4$). In \cite{WP, HPW, CLP, YP},  the authors obtained the Lagrangian densities of a $3$-uniform linear path of length $3$ or $4$, the Lagrangian densities of $\{123,234,456\}$, \{123,345,561\} and \{123,124,345\},the disjoint union of a $3$-uniform linear path of length $2$ or $3$ and a  $3$-uniform matching, and the disjoint union of a $3$-uniform tight path of length $2$  and a  $3$-uniform matching. These were all the previously known results on Lagrangian densities. For $3$-uniform graphs spanned by 3  edges, there  is one remaining unsolved case: $K_4^{3-}=\{123, 124, 134\}$.
By Proposition \ref{relationlt}, the Lagrangian density and the Tur\'an density of $K_4^{3-}$  equal. It would be very interesting to obtain the  Tur\'an density of $K_4^{3-}$ (a well-known conjecture of Frankl and F\"uredi) by determining the Lagrangian density of $K_4^{3-}$.

\subsection{ Main results and remarks}

For an $r$-graph $H$ on $t$ vertices, it is clear that $\pi_{\lambda}(H)\ge r!\lambda{(K_{t-1}^r)}$.

\begin{defi}
An $r$-graph $H$ on $t$ vertices is $\lambda$-perfect if $\pi_{\lambda}(H)= r!\lambda{(K_{t-1}^r)}$.
\end{defi}

 Theorem \ref{MStheo} implies that all $2$-graphs are $\lambda$-perfect. It is interesting to explore what kind of hypergraphs are $\lambda$-perfect. An $r$-uniform hypergraph is {\em linear} if any two edges have at most 1 vertex in common. Let $G\bigsqcup H$ denote the disjoint union of $G$ and $H$.  The following conjectures were proposed in \cite{YP}.
\begin{con} { \rm \cite{YP}}
For $r\ge 3$, there exists $n$ such that a linear $r$-graph with at least $n$ vertices is $\lambda$-perfect.
\end{con}

\begin{con}\label{con2} { \rm \cite{YP}}
For $r\ge 3$, there exists $n$ such that if $G$ and $H$ are $\lambda$-perfect $r$-graphs with at least $n$ vertices, then the disjoint union of $G$ and $H$ is $\lambda$-perfect.
\end{con}

Let $S_{2, t}$ denote the 3-graph with vertex set $\{v_1, v_2, u_1, u_2,..., u_t\}$ and edge set $\{v_1v_2u_1, v_1v_2u_2, ..., v_1v_2u_t\}$.
A result of Sidorenko in \cite{Sidorenko-89} implies that $S_{2, t}$ is $\lambda$-perfect. In \cite{YP}, we proved  that $S_{2,t}\sqcup H$ is $\lambda$-perfect if $H$ is $\lambda$-perfect and $t\geq 3$. Our first result in this paper removes  the condition that $t\geq 3$.

\begin{theo} \label{HS1}
If H is $\lambda$-perfect, then $H\sqcup S_{2, t}$ is $\lambda$-perfect for any $t\geq1$.
\end{theo}
Taking $H$ to be the 3-graph spanned by one edge or the 3-uniform linear path of length 2 or $3$ , and $t=1$,  and applying the above result, we can obtain the results in \cite{HK, JPW, HPW}.
The condition that  H is $\lambda$-perfect in the above theorem is  not necessary.

Let $F_5$ denote the 3-graph with vertex set $\{v_1, v_2, v_3, v_4, v_5\}$ and edge set $\{v_1v_2v_3, v_1v_2v_4, v_3v_4v_5\}$.
In \cite{YP}, we proved that $\pi_{\lambda}(F_5)=\frac{4}{9}$ and this implies that $F_5$ is not $\lambda$-perfect. We show that $F_5\sqcup S_{2, t}$ is $\lambda$-perfect. Let $H\sqcup\{e\}$ denote the disjoint union of $H$ and a single edge throughout the paper.

\begin{theo}\label{F5}
$F_5\sqcup\{e\}$ is $\lambda$-perfect.
\end{theo}

\begin{theo}\label{F5t}
$F_5\sqcup S_{2, t}$ is $\lambda$-perfect for $t\geq 2$.
\end{theo}

Determining the Tur\'an density of $K_4^{3-}$ (a well-known conjecture of Frankl and F\"uredi) is one of the major open problems in hypergraph Tur\'an problems. By Proposition \ref{relationlt}, the Lagrangian density and the Tur\'an density of $K_4^{3-}$  equal. It would be very interesting to obtain the  Tur\'an density of $K_4^{3-}$  by determining the Lagrangian density of $K_4^{3-}$. As a first attempt, we determine the Lagrangian density of $K_4^{3-}\sqcup e$. We show the following result.

\begin{theo}\label{K43}
$K_4^{3-}\sqcup S_{2, t}$ is $\lambda$-perfect for any $t\geq 1$.
\end{theo}

We also discuss $H\sqcup S_{2, t}$ for any $3-$graph $H$.

\begin{theo}\label{H}
Let H be a $3$-graph with  $s$ vertices. Then $H\sqcup S_{2, t}$ is $\lambda$-perfect if   $t\geq\frac{3}{2}s^2-\frac{11}{2}s+4$.
\end{theo}

Note that the property `$\lambda$-perfect' is monotone in the sense that an $r$-graph obtained by removing an edge from a  $\lambda$-perfect $r$-graph (keep the same vertex set) is $\lambda$-perfect. It's interesting to understand the relation between  the number of edges in a hypergraph and the `$\lambda$-perfect' property. We propose that the number of edges in  a hypergraph is no more than the number of edges in a linear hyperpath would guarantee the `$\lambda$-perfect' property.

\begin{con}
For $r\ge 3$, there exists $m_0$ such that for an $r$-graph $G$  with $m\ge m_0$ edges, if the number of vertices  in $G$ is at least $m(r-1)+1$, then  $G$ is $\lambda$-perfect.
\end{con}
Theorems \ref{HS1} to \ref{H} give some evidence to this conjecture. We remark that a hypergraph is not $\lambda$-perfect if it has many edges.
\begin{remark}
Let H be a 3-graph on t vertices with at least $\binom{t-1}{3}+\binom{t-2}{2}+2$ edges. Then H is not $\lambda$-perfect.
\end{remark}
\noindent{\em  Proof.} Let G be a 3-graph  with vertex set $[t]$ and  edge set $\{ijk|i, j, k\in[t-1]\}\cup\{ijt|i, j\in[t-2]\}\cup\{1(t-1)t\}$. Note that $G$ has $\binom{t-1}{3}+\binom{t-2}{2}+1$ edges , so G is H-free.  Let $x_i=\frac{1}{t-1}$ for $1\leq i\leq t-2$,  $x_j=\frac{1}{2t-2}$ for $t-1\leq j\leq t$, and  $\vec x$ be a feasible vector such that vertex i has weight $x_i$. Then $\lambda(G, \vec x)=\frac{1}{(t-1)^3}\bigg[\binom{t-1}{3}+\frac{1}{4}\bigg]>\frac{1}{(t-1)^3}\binom{t-1}{3}=\lambda(K_{t-1}^3)$. Therefore H is not $\lambda$-perfect.
\q

It might be interesting to understand whether the property `covering pairs' plays some role for the `$\lambda$-perfect' property.  Proposition \ref{relationlt} says that the property `covering pairs' is a sufficient condition for the Lagrangian density to be the same as the Tur\'an density. Is it a necessary condition?

\begin{ques}
Whether an  $r$-graph $G$ covering pairs implies that $G$ is not $\lambda$-perfect?
\end{ques}

\begin{ques}
Whether  an  $r$-graph $G$ covering pairs implies that $\pi(G)\neq \pi_{\lambda}(G)$ ?
\end{ques}

In contrast to the case $r=2$, we know very few about the set $\Pi_{r}$ for $r\ge 3$.  By Proposition \ref{relationlt}, if we get a Lagrangian density, then we can get a corresponding Tur\'an density. Sidorenko \cite{Sidorenko-89} first determined infinite many numbers in $\Pi_{r}$ for $r\ge 3$ by showing the Lagrangian densities of an infinite family of hypergraphs.  Researching on Lagrangian densities on hypergraphs might help us understand the set $\Pi_{r}$ better.

In 1999, Chung and Graham \cite{FG} proposed the conjecture that every element in $\Pi_r$ is a rational number. This conjecture was recently  disproved by  Baber and Talbot \cite{BT}, and Pikhurko \cite{Pikhurko2} independently by showing that a family of several r-graphs has irrational Tur\'an density. Baber and Talbot \cite{BT} asked whether there exists an r-graph F such that $\pi(F)$ is irrational. As mentioned earlier, the Lagrangian density and the Tur\'an density of $K_4^{3}$  equal. Determining  the Lagrangian density  of $K_4^{3}$  would lead to the solution of the  long-standing conjecture of Tur\'an on the Tur\'an density of $K_4^{3}$.  As a first attempt, we try to determine the Lagrangian  density of  $K_4^{3}\cup e$. We conjecture that the Lagrangian density of $K_4^{3}\cup e$ achieves on  the $3$-graph $S_2(n)$ with vertex set $[n]$ and  edge set $\{12i| i\in [n]\setminus\{1, 2\}\}\cup\{ijk| i\in [2], j, k\in [n]\setminus\{1, 2\}\}$ as $n\rightarrow\infty$. If this is true, then   the Lagrangian  density of  $K_4^{3}\cup e$ is $\frac{\sqrt{3}}{3}$, and this would answer the question of Baber and Talbot. But due to some difficulty, we can only show the following result and add one more irrational number to $\Pi_3$.

\begin{theo}\label{P4}
Let $\mathscr{P}=\{P_1, P_2, P_3, P_4\}$ where $P_1=\{123, 124, 134, 234, 567\}$, $P_2=\{123, 124, 134, 234, \\561, 562, 783\}$, $P_3=\{123, 124, 134, 234, 561, 562, 734\}$ and $P_4=\{123, 124, 134, 234, 561, 562, 357\}$. Then $\pi_{\lambda}(\mathscr{P})=\frac{\sqrt{3}}{3}$.
\end{theo}

As stated in Proposition \ref{relationlt}, the Lagrangian density and the Tur\'an density are the same for a large class of $r$-graphs.The advantage of transferring to the Lagrangian density is that we can assume that $G$ is dense (so covering pairs) when considering the maximum Lagrangian of an $H$-free $r$-graph $G$. This assumption makes the structural analysis `nicer' in some cases, so we hope that this method helps us  to better understand the set $\Pi_r$.

We give some properties of Lagrangians of $r$-graphs in the next section. The proofs of Theorems \ref{HS1} to \ref{H} will be given in Section \ref{sec3}. The proof of Theorem \ref{P4} will be given in Section \ref{irr}.

\section{Preliminaries}\label{sec2}

The following fact follows immediately from the definition of the Lagrangian.
\begin{fact}\label{mono}
Let $G_1$, $G_2$ be $r$-graphs and $G_1\subseteq G_2$. Then $\lambda (G_1) \le \lambda (G_2).$
\end{fact}

\begin{fact} {\em (\cite{FR})}\label{fact2}
Let $G$ be an $r$-graph on $[n]$. Let $\vec{x}=(x_1,x_2,\dots,x_n)$ be an optimum weight vector on  $G$. Then
$$ \frac{\partial \lambda (G, \vec{x})}{\partial x_i}=r\lambda(G)$$
for every $i \in [n]$ satisfying $x_i>0$.
\end{fact}

Given an $r$-graph $G$, and $i, j\in V(G),$ define $$L_G(j\setminus i)=\{e: i\notin e, e\cup\{j\}\in E(G)\:and\: e\cup\{i\}\notin E(G)\}.$$

\begin{fact}\label{symmetry}
Let $G$ be an $r$-graph on $[n]$. Let $\vec{x}=(x_1,x_2,\dots,x_n)$ be a feasible weight vector on $G$. Let $i,j\in [n]$, $i\neq j$ satisfying $L_G(i \setminus j)=L_G(j \setminus i)=\emptyset$. Let
$\vec{y}=(y_1,y_2,\dots,y_n)$ be defined by letting $y_\ell=x_\ell$ for every $\ell \in [n]\setminus \{i,j\}$ and $y_i=y_j={1 \over 2}(x_i+x_j)$.
Then $\lambda(G,\vec{y})\geq \lambda(G,\vec{x})$. Furthermore, if the pair $\{i,j\}$ is contained in an edge of $G$, $x_i>0$ for each $1\le i\le n$,  and $\lambda(G,\vec{y})=\lambda(G,\vec{x})$, then $x_i=x_j$.
\end{fact}
\noindent{\em  Proof.}
Since $L_G(i \setminus j)=L_G(j \setminus i)=\emptyset$, then
$$\lambda(G,\vec{y})-\lambda(G,\vec{x})=\sum_{\{i,j\} \subseteq e \in G}\left({(x_i+x_j)^2 \over 4}-x_ix_j\right)\prod\limits_{k\in e\setminus \{i,j\}}x_k \ge 0.$$
If the pair $\{i,j\}$ is contained in an edge of $G$ and $x_i>0$ for each $1\le i\le n$, then the equality holds only if $x_i=x_j$.
\q

An $r$-graph $G$ is {\em dense} if $\lambda (G') < \lambda (G)$ for every proper subgraph $G'$ of $G$. %This is equivalent to that all optimum weight vectors on $G$ are in the interior of ${\Delta}$, which means that no coordinate in an optimum weight vector is zero.

\begin{fact} {\em (\cite{FR})}\label{dense}
Let $G=(V,E)$ be a dense $r$-graph. Then $G$ covers pairs.
\end{fact}

Note that the converse of Fact \ref{dense} is not true.  For example, the Fano plane covers pairs but it is not dense. Indeed, many counterexamples exist by Theorem 2.1 in the paper of Talbot \cite{T}.
%\begin{fact} \label{2.6} {\em (\cite{Sidorenko-87})}
%If $G$ is a dense 3-graph on $[n]$ $(n\geq 4)$. Then $G$ contains a subgraph  isomorphic to $\{123, 124\}$.
%\end{fact}
While considering the Lagrangian density of an $r$-graph $F$, we can always reduce to consider dense $F$-free $r$-graphs.
\begin{remark} \label{remark}
Let $F$, $G$ be $r$-graphs and $G$ be $F$-free. Then (a) there exists a dense subgraph $G'$ of $G$ such that $\lambda{(G')}=\lambda{(G)}$ and $G'$ is $F$-free. (b) To show $\pi_{\lambda}(F)\leq a$, it's sufficient to show that $\lambda(G)\leq\frac{a}{r!}$ for any dense F-free r-graph. (c) To show that a t-vertex r-graph is $\lambda$-perfect, it's sufficient to show that $\lambda(G)\leq\lambda(K_{t-1}^r)$ for any dense F-free r-graph.
\end{remark}
\noindent{\em  Proof.} (a). Let $G$ be an $r$-graph on $n$ vertices. If $G$ is dense, then we are fine. If not, then we can find $G'\subset G$ such that $\lambda{(G')}=\lambda{(G)}$ and $|V(G')|<|V(G)|.$ If $G'$ is dense, then we stop. Otherwise, we continue this process until we find a dense subgraph. This process terminates since the number of vertices is reduced by at least one in each step.

(b) and (c) follows immediately from (a). \q

%%%%%%%%%%%%%%%%%%%%%%%%%%%%%%%%%%%%%%%%%%%%%%%%%%%%%

\section{The disjoint union of a 3-graph and $S_{2, t}$}\label{sec3}

We will prove Theorems \ref{HS1} to \ref{H} in this section.

\subsection{Preliminaries}

For a 3-graph G and $v\in V(G)$, we define the link graph of v as $G_v=\{ab|vab\in E(G)\}$. Let $\omega(G)$ denote the order of a maximum clique of G.

\begin{claim}\label{3.1}
Let G be a 3-graph with $\lambda(G)>\lambda(K_{k+1}^3)$ and let $\vec{x}$ be an optimal weight vector. Then for any $v\in V(G)$, its weight $x_v$ satisfies that $x_v<1-\frac{\sqrt{k(k-1)}}{k+1}$.
\end{claim}
\emph{Proof.} Since the link graph $G_v$ of v is a graph, by Fact \ref{fact2} and the theorem of Motzkin and Straus (Theorem \ref{MStheo}), we have $$3\lambda(K_{k+1}^3)<3\lambda(G)=\frac{\partial\lambda(G)}{\partial x_v}\leq\frac{1}{2}(1-x_v)^2.$$ So $$x_v<1-\frac{\sqrt{k(k-1)}}{k+1}.$$
\q

\begin{claim}\label{3.2}
Let G be a 3-graph with $\lambda(G)>\lambda(K_{k+1}^3)$ and let $\vec{x}$ be an optimal weight vector. Then for $v\in V(G)$ with $\omega(G_v)\leq k$, its weight $x_v$ satisfies that $x_v<\frac{1}{k+1}$.
\end{claim}
\emph{Proof.} Since the link graph $G_v$ of v is a graph, by Fact \ref{fact2} and the theorem of Motzkin and Straus (Theorem \ref{MStheo}), we have $$3\lambda(K_{k+1}^3)<3\lambda(G)=\frac{\partial\lambda(G)}{\partial x_v}\leq\frac{1}{2}\times(1-\frac{1}{k})(1-x_v)^2.$$ So$$x_v<\frac{1}{k+1}.$$
\q

\begin{claim}\label{3.3}
Let v be a vertex in a 3-graph G and $x_v$ be the weight of v in an optimal weight vector $\vec x$ of G. If $G-\{v\}$ is H-free, then $\lambda(G)\leq\frac{\pi_{\lambda}(H)(1-x_v)^3}{6(1-3x_v)}$.
\end{claim}
\emph{Proof.} Since $G-\{v\}$ is H-free, then $\lambda(G-\{v\}, \vec x)\leq\frac{\pi_{\lambda}(H)}{6}(1-x_v)^3$. Therefore $$\lambda(G)\leq\frac{\pi_{\lambda}(H)}{6}(1-x_v)^3+x_v\frac{\partial{\lambda(G)}}{\partial x_v}.$$ By Fact \ref{fact2}, we have $$\lambda(G)\leq\frac{\pi_{\lambda}(H)}{6}(1-x_v)^3+3x_v\lambda(G).$$ Then $\lambda(G)\leq\frac{\pi_{\lambda}(H)(1-x_v)^3}{6(1-3x_v)}$.
\q

\begin{remark} \label{remark33}
$f(x)=\frac{(1-x)^3}{1-3x}$ is increasing in $(0, \frac{1}{3}).$
\end{remark}
\emph{Proof.} Since $f'(x)=\frac{6x(1-x)^2}{(1-3x)^2}$, then $f'(x)>0$ in $(0, \frac{1}{3}).$ So $f(x)=\frac{(1-x)^3}{1-3x}$ is increasing in $(0, \frac{1}{3}).$
\q

\begin{claim}\label{3.4}
Let a 3-graph G be $H\cup S_{2, t}$-free, where H is a 3-graph with s vertices. Let $S_{2, s+t}=\{v_1v_2b_1, v_1v_2b_2, ..., v_1v_2b_{s+t} \}\subseteq G$. Then $G-\{v_1, v_2\}$ is H-free.
\end{claim}
\emph{Proof.} Suppose that $H\subseteq G-\{v_1, v_2\}$. Since $|V(H)|=s$, then $|\{b_1, b_2, ..., b_{s+t}\}\cap V(H)|\leq s$, and $|\{b_1, ..., b_{s+t}\}\setminus V(H)|\geq t$. So the induced subgraph of G by $\{v_1, v_2, b_1, ... , b_{s+t}\}\setminus V(H)$ contains an $S_{2, t}$, and G contains $H\cup S_{2, t}.$
\q

\begin{claim}\label{3.5}
Let a 3-graph G be $H\cup S_{2, t}$-free, where H is a 3-graph with s vertices. Let $v\in V(H)$. If $H\subseteq G-\{v\}$, then $\omega(G_{v})\leq s+t$.
\end{claim}
\emph{Proof.} Since $G-\{v\}$ contains H. If $\omega(G_{v})\geq s+t+1,$ and assume that a maximum clique of $G_{v}$ has vertex set $U_1$. Since $|U_1\cap V(H)|\leq s,$ then $|U_1\setminus V(H)|\geq t+1$ and the induced subgraph of G by $U_1\cup\{v_1\}\setminus V(H)$ contains an $S_{2, t}$. So $H\cup S_{2, t}\subseteq G$, a contradiction.
\q

\begin{claim}\label{3.6}
Let a 3-graph G be $H\cup S_{2, t}$-free, where H is a 3-graph with s vertices. If $H\subseteq G-\{v_1\}$ and $H\nsubseteq G-\{v_1, v_2\}$, then $\omega((G-\{v_2\})_{v_1})\leq s+t-1.$
\end{claim}
\emph{Proof.} Assume that $\omega((G-\{v_2\})_{v_1})\geq s+t$ and a maximum clique of $(G-\{v_2\})_{v_1}$ has vertex set $U_2$. Since $H\subseteq G-\{v_1\}$ and $H\nsubseteq G-\{v_1, v_2\}$, then $v_2\in V(H)$. Therefore $|U_2\cap V(H)|\leq s-1.$ So $|U_2\setminus V(H)|\geq t+1$ and the induced subgraph of G by $U_2\cup\{v_1\}\setminus V(H)$ contains an $S_{2, t}$. Therefore $H\cup S_{2, t}\subseteq G$.
\q

\begin{theo}{\em (\cite{BT1}\cite{Raz})}\label{FL}
$\pi(K_4^{3-})\leq 0.2871$, $\pi(K_4^3)\leq0.5617.$
\end{theo}

\subsection{The disjoint union of a $\lambda$-perfect 3-graph and $S_{2, t}$}

We show that the disjoint union of a $\lambda$-perfect 3-graph and $S_{2, t}$ is $\lambda$-perfect.

\emph{Proof of Theorem \ref{HS1}.} Assume that H is $\lambda$-perfect on $s\geq3$ vertices. Note that $H\cup S_{2, t}$ has s+t+2 vertices. By Remark \ref{remark}(c), it's sufficient to show that if G is $H\sqcup S_{2, t}$-free dense 3-graph then  $\lambda(G)\leq\lambda(K_{s+t+1}^3)$. Suppose on the contrary that $\lambda(G)>\lambda(K_{s+t+1}^3)$. Let $\vec x$ be an optimal weight vector of G.

\emph{Case 1.} There exists $v\in V(G)$ with weight $x_v$ such that $G-\{v\}$ is H-free.

By Claim \ref{3.1}, $x_v<1-\frac{\sqrt{(s+t)(s+t-1)}}{s+t+1}$. By Claim \ref{3.3} and that H is $\lambda$-perfect, $$\lambda(G)\leq\frac{\lambda(K_{s-1}^3)(1-x_v)^3}{1-3x_v}=f(x_v).\eqno {(3)}.$$
By Remark \ref{remark33}, $f(x_v)$ is increasing in $[0, 1-\frac{\sqrt{(s+t)(s+t-1)}}{s+t+1})$, then
\begin{eqnarray*}
\lambda(G)&\leq&f(1-\frac{\sqrt{(s+t)(s+t-1)}}{s+t+1})\\
&=&\frac{1}{6}\frac{(s-2)(s-3)}{(s-1)^2}\frac{\frac{(s+t)(s+t-1)\sqrt{(s+t)(s+t-1)}}{(s+t+1)^3}}{\frac{3\sqrt{(s+t)(s+t-1)}}{s+t+1}-2}\\
&=&\frac{1}{6}\frac{(s+t)(s+t-1)}{(s+t+1)^2}\times\frac{(s-2)(s-3)}{(s-1)^2}\times\frac{\sqrt{(s+t)(s+t-1)}}{3\sqrt{(s+t)(s+t-1)}-2(s+t+1)}\\
&=&\lambda(K_{s+t+1}^3)\times\frac{(s-2)(s-3)}{(s-1)^2}\times\frac{\sqrt{(s+t)(s+t-1)}}{3\sqrt{(s+t)(s+t-1)}-2(s+t+1)}.
\end{eqnarray*}
To prove $\lambda(G)\leq \lambda(K_{s+t+1}^3)$, it's sufficient to prove that $$\frac{(s-2)(s-3)}{(s-1)^2}\times\frac{\sqrt{(s+t)(s+t-1)}}{3\sqrt{(s+t)(s+t-1)}-2(s+t+1)}\leq 1.$$
This is equivalent to $$\frac{2(s-1)^2(s+t+1)}{2s^2-s-3}\leq\sqrt{(s+t)(s+t-1)}.$$
The above inequality is equivalent to $$s+t-\frac{1}{2}-\frac{(6t-1)s-10t-1}{4s^2-2s-6}\leq\sqrt{(s+t)(s+t-1)}.$$
By direct calculation, it holds for $s$=3 or 4 and $t=1$.
Since $$s+t-\frac{1}{2}-\frac{(6t-1)s-10t-1}{4s^2-2s-6}\leq s+t-\frac{1}{2}-\frac{1}{s+1},$$ and $$s+t-\frac{1}{2}-\frac{1}{s+1}\leq\sqrt{(s+t)(s+t-1)}$$ holds for any $s\geq 3$ and $t\geq2$ or $s\geq 5$ and $t=1$. So $\lambda(G)\leq \lambda(K_{s+t+1}^3)$, a contradiction.

\emph{Case 2.} For any $v\in V(G)$, $H\subseteq G-\{v\}$.

Since $\lambda(G)>\lambda(K_{s+t+1}^3)$ and $S_{2, s+t}$ is $\lambda$-perfect, then $S_{2, s+t}=\{v_1v_2b_1, v_1v_2b_2, ..., v_1v_2b_{s+t}\}\subseteq G$. By Claim \ref{3.4}, $G-\{v_1, v_2\}$ is H-free. By Claim \ref{3.5}, we have $$\omega(G_{v_1})\leq s+t\quad{\rm and}\quad\omega(G_{v_2})\leq s+t.$$ By Claim \ref{3.6}, we have $\omega((G-\{v_2\})_{v_1})\leq s+t-1$ and $\omega((G-\{v_1\})_{v_2})\leq s+t-1$.

Assume the weight of $v_1$ and $v_2$ are $a_1$ and $a_2$ respectively, and $a_1+a_2=2a$. Since $G-\{v_1, v_2\}$ is H-free and H is $\lambda$-perfect, the contribution of edges containing neither $v_1$ nor $v_2$ to $\lambda(G, \vec x)$ is at most $\lambda(K_{s-1}^3)(1-2a)^3$. Since $\omega((G-\{v_2\})_{v_1})\leq s+t-1$ and $\omega((G-\{v_1\})_{v_2})\leq s+t-1$, by Theorem \ref{MStheo}, the contribution of edges containing either $v_1$ or $v_2$ to $\lambda(G, \vec x)$ is at most $2\times\frac{1}{2}a(1-\frac{1}{s+t-1})(1-2a)^2$. The contribution of edges containing both $v_1$ and $v_2$ to $\lambda(G, \vec x)$ is at most $a^2(1-2a)$. Therefore
\begin{eqnarray*}
\lambda(G)&\leq&\lambda(K_{s-1}^3)(1-2a)^3+a^2(1-2a)+2\times\frac{1}{2}a(1-\frac{1}{s+t-1})(1-2a)^2\\%\quad (a\leq\frac{1}{t+2})
&\leq&\lambda(K_{s+t-1}^3)(1-2a)^3+a^2(1-2a)+2\times\frac{1}{2}a(1-\frac{1}{s+t-1})(1-2a)^2\\
&=&\lambda(K_{s+t+1}^3, (a, a, \frac{1-2a}{s+t-1}, ..., \frac{1-2a}{s+t-1}))\\
&\leq&\lambda(K_{s+t+1}^3).
\end{eqnarray*}
\q

\subsection{Disjoint union of $F_5$ and $S_{2, t}$}

We show that the disjoint union of $F_5$ and $S_{2, t}$ is $\lambda$-perfect.

\emph{Proof of Theorem \ref{F5t}.} Let $t\geq 2.$ Let G be a dense $F_5\sqcup S_{2, t}$-free 3-graph on n vertices. By Remark \ref{remark}(c), it's sufficient to show that $\lambda(G)\leq \lambda(K_{t+6}^3)$. Suppose on the contrary that $\lambda(G)>\lambda(K_{t+6}^3)$. By Claim \ref{3.1}, for any $v\in V(G)$ with weight $x_v$ we have $$x_v<1-\frac{\sqrt{(t+5)(t+4)}}{t+6}.$$

\emph{Case 1.} $G-\{v\}$ is $F_5$-free for some $v\in V(G)$.

By the result in \cite{YP} that $\pi_{\lambda}(F_5)=\frac{4}{9}$ and Claim \ref{3.3}, we have
\begin{eqnarray*}
\lambda(G)&\leq&\frac{2}{27}\frac{(1-x_v)^3}{1-3x_v}=f(x_v)\\
&\leq&f(1-\frac{\sqrt{(t+5)(t+4)}}{t+6})\\
&=&\frac{1}{6}\frac{(t+4)(t+5)}{(t+6)^2}\times\frac{4}{9}\times\frac{\sqrt{(t+5)(t+4)}}{3\sqrt{(t+5)(t+4)}-2(t+6)}\\
&=&\lambda(K_{t+6}^3)\times\frac{4}{9}\times\frac{\sqrt{(t+5)(t+4)}}{3\sqrt{(t+5)(t+4)}-2(t+6)}.
\end{eqnarray*}
It's sufficient to show that $$\frac{4}{9}\times\frac{\sqrt{(t+5)(t+4)}}{3\sqrt{(t+5)(t+4)}-2(t+6)}\leq 1.$$
It's equivalent to show that $$18(t+6)\leq23\sqrt{(t+5)(t+4)}.$$ It holds for $t\geq 2.$

\emph{Case 2.} $G-\{v\}$ contains $F_5$ for any $v\in V(G)$.

Since $\lambda(G)>\lambda(K_{t+6}^3)$ and $S_{2, t+5}$ is $\lambda$-perfect, then $S_{2, t+5}=\{v_1v_2u_1, v_1v_2u_2, ..., v_1v_2u_{t+5}\}\subseteq G$.  Assume the weight of $v_1$ and $v_2$ are $a_1$ and $a_2$ respectively, and $a_1+a_2=2a$. By Claim \ref{3.4}, $G-\{v_1, v_2\}$ is $F_5$-free. So the contribution of edges containing neither $v_1$ nor $v_2$ to $\lambda(G, \vec x)$ is at most $\frac{\pi_{\lambda}(F_5)}{6}(1-2a)^3=\frac{2}{27}(1-2a)^3$. By Claim \ref{3.6}, we have $\omega((G-\{v_2\})_{v_1})\leq t+4$ and $\omega((G-\{v_1\})_{v_2})\leq t+4$. By Theorem \ref{MStheo}, the contribution of edges containing either $v_1$ or $v_2$ to $\lambda(G, \vec x)$ is at most $2a\times\frac{1}{2}(1-\frac{1}{t+4})(1-2a)^2$. The contribution of edges containing both $v_1$ and $v_2$ to $\lambda(G, \vec x)$ is at most $a^2(1-2a)$. Therefore
\begin{eqnarray*}
\lambda(G)&\leq&\frac{2}{27}(1-2a)^3+a^2(1-2a)+2a\times\frac{1}{2}(1-\frac{1}{t+4})(1-2a)^2\\
&\leq&\lambda(K_{t+4}^3)(1-2a)^3+a^2(1-2a)+2a\times\frac{1}{2}(1-\frac{1}{t+4})(1-2a)^2\\
&=&\lambda(K_{t+6}^3, (a, a, \frac{1-2a}{t+4}, ..., \frac{1-2a}{t+4}))\\
&\leq&\lambda(K_{t+6}^3).
\end{eqnarray*}
\q
\bigskip

The following lemma is prepared for the proof of Theorem \ref{F5}.

\begin{lemma}\label{1}
Let G be a dense and $F_5\sqcup\{e\}$-free 3-graph. If $K_4^{3-}\sqcup K_4^{3-}\subseteq G$, then $\lambda(G)<\frac{5}{49}=\lambda(K_7^3).$
\end{lemma}
\emph{Proof.} Let $K=K_4^{3-}\sqcup K_4^{3-}$ with vertex set $\{a_1, a_2, b_1, b_2, ..., b_6\}$ and edge set $\{a_1b_1b_2, a_1b_1b_3, a_1b_2b_3, a_2b_4b_5,\\ a_2b_4b_6, a_2b_5b_6\}$. Let $\vec x$ be an optimal weight vector. To simplify the notation, we simply write the weight of vertex $a$ in $\vec x$ as $a$. In other words, in the proof, $a$ sometimes means vertex $a$, sometimes means the weight of vertex $a$ in the optimal weight vector $\vec x$. By Fact \ref{fact2}, $$24\lambda=\sum_{i=1}^2\frac{\partial\lambda}{\partial a_i}+\sum_{i=1}^6\frac{\partial\lambda}{\partial b_i}. \eqno (1)$$

\begin{claim}
i) For vertices $c_i, c_j\in V(G)\setminus V(K)$, the product $c_ic_j$ appears at most twice in $\sum_{i=1}^2\frac{\partial\lambda}{\partial a_i}+\sum_{i=1}^6\frac{\partial\lambda}{\partial b_i}$.

ii) For vertices $c_i \in V(G)\setminus V(K)$, the product $c_ib_j$ appears at most twice in $\sum_{i=1}^2\frac{\partial\lambda}{\partial a_i}+\sum_{i=1}^6\frac{\partial\lambda}{\partial b_i}$.

iii) For vertices $c_i \in V(G)\setminus V(K)$, the product $c_ia_1$, $c_ia_2$ appear at most 7 times in $\sum_{i=1}^2\frac{\partial\lambda}{\partial a_i}+\sum_{i=1}^6\frac{\partial\lambda}{\partial b_i}$.

iv) The product $b_ib_j$ appears at most 3 times in $\sum_{i=1}^2\frac{\partial\lambda}{\partial a_i}+\sum_{i=1}^6\frac{\partial\lambda}{\partial b_i}$.

v) The product $a_ib_j$, $a_1a_2$ appear at most 6 times in $\sum_{i=1}^2\frac{\partial\lambda}{\partial a_i}+\sum_{i=1}^6\frac{\partial\lambda}{\partial b_i}$.
\end{claim}
\emph{Proof.}  i)Let F be a copy of $F_5\cup\{e\}$. We first show that $c_ic_j$ can appear at most 1 time in $\frac{\partial\lambda}{\partial a_1}+\sum_{i=1}^3\frac{\partial\lambda}{\partial b_i}$. If not, suppose that $c_ic_j$ appears at least 2 times in it. If $c_ic_jb_s, c_ic_jb_t\in E(G)$ where $b_s, b_t\in \{b_1, b_2, b_3\}$, then $\{a_1, b_s, b_t, c_i, c_j\}$ forms an $F_5$. So there is an $F$ formed by $F_5\sqcup\{a_2b_4b_5\}$. If $c_ic_jb_s, c_ic_ja_1\in E(G)$ where $b_s\in \{b_1, b_2, b_3\}$, then $\{a_1, b_t, b_s, c_i, c_j\}$ forms an $F_5$, where $b_t\in \{b_1, b_2, b_3\}\setminus\{b_s\}$. So $F_5\sqcup\{a_2b_4b_5\}\subseteq G$. Similarly, $c_ic_j$ can appear at most 1 time in $\frac{\partial\lambda}{\partial a_2}+\sum_{i=4}^6\frac{\partial\lambda}{\partial b_i}$. Therefore $c_ic_j$ appears at most twice in $\sum_{i=1}^2\frac{\partial\lambda}{\partial a_i}+\sum_{i=1}^6\frac{\partial\lambda}{\partial b_i}$.\\

ii) Without loss of generality, let $b_j=b_1$. We first show that $c_ib_1$ appears in $\frac{\partial\lambda}{\partial a_1}+\frac{\partial\lambda}{\partial b_2}+\frac{\partial\lambda}{\partial b_3}$ at most 1 time. Since if $c_ib_1b_2, c_ib_1b_3\in E(G)$, then $\{a_1, b_1, b_2, b_3, c_i\}$ forms an $F_5$, then there is an $F$ formed by $F_5\sqcup\{a_2b_4b_5\}$. If $c_ib_1b_2$ (or $c_ib_1b_3$) and $c_ib_1a_1\in E(G)$, then $\{a_1, b_1, b_2, b_3, c_i\}$ forms an $F_5$, then there is an $F$ formed by $F_5\sqcup\{a_2b_4b_5\}$.

Next we show that $c_ib_1$ appears in $\frac{\partial\lambda}{\partial a_2}+\sum_{i=4}^6\frac{\partial\lambda}{\partial b_i}$ at most 1 time. If $c_ib_1b_s, c_ib_1b_t\in E(G)$, where $b_s, b_t\in\{b_4, b_5, b_6\}$, then $\{a_2, b_1, b_s, b_t, c_i\}$ forms an $F_5$, and there is an $F$ formed by $F_5\sqcup\{a_1b_2b_3\}$. If $c_ib_1a_2, c_ib_1b_s\in E(G)$, where $b_s\in\{b_4, b_5, b_6\}$, then $\{a_2, b_1, b_s, b_t, c_i\}$ forms an $F_5$ where $b_s, b_t\in\{b_4, b_5, b_6\}$. So there is an $F$ formed by $F_5\sqcup\{a_1b_2b_3\}$.\\

iii) Let $j=1$ or 2. Since $c_ia_j\notin \frac{\partial\lambda}{\partial a_j}$, then $c_ia_j$ appears at most 7 times.\\

iv) If $b_i\in \{b_1, b_2, b_3\}$ and $b_j\in \{b_4, b_5, b_6\}$, without loss of generality, let $i=1$ and $j=4$. Note that $b_1b_2b_4\notin E(G)$. Otherwise $\{a_1, b_1, b_2, b_3, b_4\}$ forms an $F_5$ and there is an $F$ formed by $F_5\sqcup\{a_2b_5b_6\}$. Similarly, $b_1b_3b_4, b_1b_4b_5, b_1b_4b_6, b_1b_2b_5, b_1b_2b_6\notin E(G)$. So at most $b_1b_4a_1, b_1b_4a_2\in E(G)$. If $b_i, b_j\in \{b_1, b_2, b_3\}$, without loss of generality, let $i=1$ and $j=2$. Then at most $b_1b_2a_1, b_1b_2a_2, b_1b_2b_3\in E(G)$.\\

v) Since $a_ib_j$ can not appear in $\frac{\partial\lambda}{\partial a_j}+\frac{\partial\lambda}{\partial b_i}$, then $a_ib_j$ can appear at most 6 times in $\sum_{i=1}^2\frac{\partial\lambda}{\partial a_i}+\sum_{i=1}^6\frac{\partial\lambda}{\partial b_i}$. The reason for $a_1a_2$ is similar.
\q

Combining (1) and Claim \ref{3.3}, we have $$24\lambda\leq2\sum_{c_i,c_j\in V(G)\setminus V(K)}c_ic_j+2\sum_{i\in [6], c_j\in V(G)\setminus V(K)}b_ic_j+7\sum_{i\in [2], c_j\in V(G)\setminus V(K)}a_ic_j+$$ $$3\sum_{i\neq j\in[6]}b_ib_j+\\6\sum_{i=1}^2\sum_{j=1}^6a_ib_j+6a_1a_2.$$
Let $\sum_{i}c_i=c$, $\sum_{i=1}^2a_i=a$, $\sum_{i=1}^6b_i=b$, then $a+b+c=1$. So
\begin{eqnarray*}
24\lambda&\leq &2\bigg(\frac{c}{n-8}\bigg)^2\binom{n-8}{2}+2bc+7ac+3\bigg(\frac{b}{6}\bigg)^2\binom{6}{2}+6ab+6\bigg(\frac{a}{2}\bigg)^2\\
&\leq&c^2+2bc+7ac+\frac{3}{2}b^2+6ab+\frac{3}{2}a^2\\
&=&(1-a-b)^2+2b(1-a-b)+7a(1-a-b)+\frac{3}{2}b^2+6ab+\frac{3}{2}a^2\\
&=&-\frac{9}{2}a^2+\frac{1}{2}b^2-ab+5a+1\\
&=&-\frac{9}{2}a^2+(5-b)a+\frac{1}{2}b^2+1\\
&\leq&\frac{5}{9}b^2-\frac{5}{9}b+\frac{43}{18}\qquad(a=\frac{5-b}{9})\\
&\leq&\frac{43}{18}.\\
\end{eqnarray*}
So
$$\lambda(G)\leq\frac{43}{18\times 24}<\frac{5}{49}=\lambda(K_7^3).$$
\q
\\
\emph{Proof of Theorem \ref{F5}.} Let G be a dense $F_5\cup\{e\}$-free 3-graph. By Remark \ref{remark}(c), it's sufficient to show that $\lambda(G)\leq \lambda(K_7^3)$. Suppose on the contrary that $\lambda(G)>\frac{5}{49}=\lambda(K_7^3)$. For $v\in V(G)$ with weight $x_v$, if $G-\{v\}$ is $K_4^{3-}$-free, then by Proposition \ref{relationlt}, Theorem \ref{3.3} and Claim \ref{FL}, we have $$\lambda(G)\leq\frac{\pi_{\lambda}(K_4^{3-})}{6}\frac{(1-x_v)^3}{1-3x_v}=\frac{\pi(K_4^{3-})}{6}\frac{(1-x_v)^3}{1-3x_v}\leq\frac{1}{18}\frac{(1-x_v)^3}{1-3x_v}.$$ By Claim \ref{3.1}, we have $x_v<1-\frac{\sqrt{30}}{7}<\frac{1}{3}.$ By Remark \ref{remark33}, $\frac{1}{18}\frac{(1-x_v)^3}{1-3x_v}$ is increasing in $(0, \frac{1}{3})$, so $$\lambda(G)\leq\frac{1}{18}\frac{(1-x_v)^3}{1-3x_v}\bigg|_{x_v=1-\frac{\sqrt{30}}{7}}<\frac{5}{49}.$$
Hence we may assume that $G-\{v\}$ still contains $K_4^{3-}$ for any vertex $v$. Note that $\omega(G_v)\leq6$. Since if $\omega(G_v)\geq 7$, assume that $U=\{u_1, u_2, ..., u_7\}$ is a clique in $G_v$. Let $K_4^{3-}\subseteq G-\{v\}$ have the vertex set W. Then $|U-W|\geq 3$ and we can find a $K_4^{3-}$ in $\{v_1\}\cup U-W$. Therefore $K_4^{3-}\sqcup K_4^{3-}\subseteq G$. So by Fact \ref{fact2} and the theorem of Motzkin and Straus (Theorem \ref{MStheo}), we have $$3\times\frac{5}{49}<3\lambda=\frac{\partial\lambda}{\partial x_v}\leq\frac{1}{2}(1-\frac{1}{6})(1-x_v)^2.$$ So $x_v<\frac{1}{7}$ holds for any $v$.

If for $v\in V(G)$ with weight $x_v$, $G-\{v\}$ is $F_5$-free, then by the result in \cite{YP} that $\pi_{\lambda}(F_5)=\frac{4}{9}$ and Claim \ref{3.3}, we have $$\lambda(G)\leq\frac{2}{27}\frac{(1-x_v)^3}{1-3x_v}\leq\frac{2}{27}\frac{(1-x_v)^3}{1-3x_v}\bigg|_{x_v=\frac{1}{7}}<\frac{5}{49}.$$
So we may assume that $G-\{v\}$ still contains both $F_5$ and $K_4^{3-}$ for any vertex $v$. Since $\lambda(G)>\frac{5}{49}$ and $S_{2, t}$ is $\lambda$-perfect, then $S_{2, 6}\subseteq G.$ Let $S_{2, 6}=\{v_1v_2b_1, v_1v_2b_2, ..., v_1v_2b_6\}$. By Claim \ref{3.4}, $G-\{v_1, v_2\}$ is $F_5$-free. By Claim \ref{3.6}, $\omega((G-\{v_2\})_{v_1})\leq 5$ and $\omega((G-\{v_1\})_{v_2})\leq 5$. Assume the weight of $v_1$ and $v_2$ are $a_1$ and $a_2$ respectively and $a_1+a_2=2a$. By Claim \ref{3.2}, $a_1<\frac{1}{7}$ and $a_1<\frac{1}{7}$. Therefore $a<\frac{1}{7}$. So
\begin{eqnarray*}
\lambda(G)&\leq &a^2(1-2a)+\frac{2}{27}(1-2a)^3+2a\times\frac{1}{2}\times(1-\frac{1}{5})(1-2a)^2\\
&=&\frac{82}{135}a^3-\frac{59}{45}a^2+\frac{16}{45}a+\frac{2}{27}=f(a)\\
f'(a)&=&\frac{82}{45}a^2-\frac{118}{45}a+\frac{16}{45}>0\\
\end{eqnarray*}
if $a\in [0, \frac{1}{7}]$. So $f(a)$ is increasing in $[0, \frac{1}{7}]$. Then $$\lambda(G)\leq f(a)=f(\frac{1}{7})<\frac{5}{49}.$$
\q

%%\begin{coro}
%%$\pi_{\lambda}(F=F_5\sqcup\{e\})=\frac{30}{49}$ i.e. F is a $\lambda$-perfect 3-graph.
%%\end{coro}
%%\emph{Proof.} Since $K_7^3$ is F-free with $\lambda(K_7^3)=\frac{5}{49}$, So $\pi_{\lambda}(F)=3!\times\frac{5}{49}=\frac{30}{49}.$
%%\q
\subsection{The disjoint union of $K_4^{3-}$ and $S_{2, t}$}

We show that the disjoint union of $K_4^{3-}$ and $S_{2, t}$ is $\lambda$-perfect.

\emph{Proof of Theorem \ref{K43}.} Let G be a dense $K_4^{3-}\sqcup S_{2, t}$-free 3-graph on n vertices. By Remark \ref{remark}(c), it's sufficient to show that $\lambda(G)\leq \lambda(K_{t+5}^3)$. Suppose on the contrary that $\lambda(G)>\lambda(K_{t+5}^3)$.
If there exists $v\in V(G)$ with weight $x_v$ such that $G-\{v\}$ is $S_{2, t}$-free, then by Claim \ref{3.1}, Claim \ref{3.3} and $S_{2, t}$ is $\lambda$-perfect, we have $$x_v<1-\frac{\sqrt{(t+3)(t+4)}}{t+5}\quad {\rm and} \quad\lambda(G)\leq\frac{\pi_{\lambda}(S_{2, t})(1-x_v)^3}{6(1-3x_v)}\leq\frac{\lambda(K_{t+1}^3)(1-x_v)^3}{1-3x_v}=f(x_v).$$
Since $f(x)$ is increasing in $(0, \frac{1}{3})$, then
\begin{eqnarray*}
\lambda(G)
&\leq&f(1-\frac{\sqrt{(t+3)(t+4)}}{t+5})\\
&=&\frac{1}{6}\frac{t(t-1)}{(t+1)^2}\frac{\frac{(t+4)(t+3)\sqrt{(t+4)(t+3)}}{(t+5)^3}}{\frac{3\sqrt{(t+4)(t+3)}}{t+5}-2}\\
&=&\frac{1}{6}\frac{(t+4)(t+3)}{(t+5)^2}\frac{t(t-1)}{(t+1)^2}\frac{\sqrt{(t+4)(t+3)}}{3\sqrt{(t+4)(t+3)}-2(t+5)}\\
&=&\lambda(K_{t+5}^3)\frac{t(t-1)}{(t+1)^2}\frac{\sqrt{(t+4)(t+3)}}{3\sqrt{(t+4)(t+3)}-2(t+5)}.
\end{eqnarray*}
So it's sufficient to show that $$\frac{t(t-1)}{(t+1)^2}\times\frac{\sqrt{(t+4)(t+3)}}{3\sqrt{(t+4)(t+3)}-2(t+5)}\leq 1.$$
This is equivalent to show that $$t+\frac{7}{2}-\frac{11t+1}{4t^2+14t+6}\leq\sqrt{(t+4)(t+3)}.$$
This is true since $$\frac{11t+1}{4t^2+14t+6}\geq\frac{1}{t+1}, (t\geq 1)$$ and $$t+\frac{7}{2}-\sqrt{(t+4)(t+3)}=\frac{(t+\frac{7}{2})^2-(t+3)(t+4)}{t+\frac{7}{2}+\sqrt{(t+4)(t+3)}}<\frac{1}{4(t+1)}<\frac{1}{t+1}.$$

So we may assume that $S_{2, t}\in G-\{v\}$ for any $v\in V(G)$. By Claim \ref{3.5}, $\omega(G_v)\leq t+4$. Let $x_v$ be the weight of v in an optimum vector. By Claim \ref{3.2}, $x_v<\frac{1}{t+5}$.

\emph{Case 1.} $G-\{v\}$ is $K_4^{3-}$-free for some $v\in V(G)$.

Then by Proposition \ref{relationlt}, Claim \ref{3.3} and Theorem \ref{FL}, we have
\begin{eqnarray*}
\lambda(G)&\leq&\frac{1}{18}\frac{(1-x_v)^3}{1-3x_v}=f(x_v) \quad (x_v<\frac{1}{t+5})\\
&<&f(\frac{1}{t+5})\\
&=&\frac{1}{18}\frac{(t+4)^3}{(t+2)(t+5)^2}\\
&<&\frac{1}{6}\frac{(t+4)(t+3)}{(t+5)^2}=\lambda(K_{t+5}^3).
\end{eqnarray*}

\emph{Case 2.} $G-\{v\}$ contains $K_4^{3-}$ for any $v\in V(G)$.

Since $\lambda(G)>\lambda(K_{t+5}^3)$ and $S_{2, t+4}$ is $\lambda$-perfect, then $S_{2, t+4}=\{v_1v_2u_1, v_1v_2u_2, ..., v_1v_2u_{t+4}\}\subseteq G$. By Claim \ref{3.4}, $G-\{v_1, v_2\}$ is $K_4^{3-}$-free. Since $G-\{v_1\}$ contains $K_4^{3-}$ and $S_{2, t}$, by Claim \ref{3.6}, $\omega((G-\{v_2\})_{v_1})\leq t+3$ and $\omega((G-\{v_1\})_{v_2})\leq t+3$. Assume the weight of $v_1$ and $v_2$ are $a_1$ and $a_2$ respectively. Let $a_1+a_2=2a$. So
\begin{eqnarray*}
\lambda(G)&\leq &\frac{1}{18}(1-2a)^3+a^2(1-2a)+2a\times\frac{1}{2}(1-\frac{1}{t+3})(1-2a)^2\\
&\leq &\lambda(K_{t+3}^3)(1-2a)^3+a^2(1-2a)+2a\times\frac{1}{2}(1-\frac{1}{t+3})(1-2a)^2\\
&<&\lambda(K_{t+5}^3, (a, a, \frac{1-2a}{t+3}, ..., \frac{1-2a}{t+3}))\\
&\leq&\lambda(K_{t+5}^3).
\end{eqnarray*}
\q

\subsection{General case}

We show that $H\sqcup S_{2, t}$ is $\lambda$-perfect if H is a $3$-graph with  $s$ vertices and $t\geq\frac{3}{2}s^2-\frac{11}{2}s+4$.

\emph{Proof of Theorem \ref{H}.} Let G be a dense $H\sqcup S_{2, t}$-free 3-graph on n vertices.  By Remark \ref{remark}(c), it's sufficient to show that $\lambda(G)\leq\lambda(K_{s+t+1}^3)$. Suppose on the contrary that $\lambda(G)>\lambda(K_{s+t+1}^3)$.
If there exists $v\in V(G)$ with weight $x_v$ such that $G-\{v\}$ is $H$-free, then by Claim \ref{3.1} and Claim \ref{3.3}, we have $$\lambda(G)\leq\frac{1}{6}\pi_{\lambda}(H)\frac{(1-x_v)^3}{1-3x_v},\quad{\rm and}\quad x_v<1-\frac{\sqrt{(s+t)(s+t-1)}}{s+t+1}.$$
Since $\pi_{\lambda}(H)\leq\pi_{\lambda}(K_s^3)\leq 1-\frac{2}{(s-1)(s-2)}$ (see \cite{Keevash}), then
\begin{eqnarray*}
\lambda(G)&\leq&\frac{1}{6}\bigg(1-\frac{2}{(s-1)(s-2)}\bigg)\frac{(1-x_v)^3}{1-3x_v}\bigg|_{x_v=1-\frac{\sqrt{(s+t)(s+t-1)}}{s+t+1}}\\
&\leq&\lambda(K_{s+t+1}^3)\frac{s^2-3s}{(s-1)(s-2)}\times\frac{\sqrt{(s+t)(s+t-1)}}{3\sqrt{(s+t)(s+t-1)}-2}.
\end{eqnarray*}
So it's sufficient to show that $$\frac{s^2-3s}{(s-1)(s-2)}\times\frac{\sqrt{(s+t)(s+t-1)}}{3\sqrt{(s+t)(s+t-1)}-2}\leq 1.$$
This is equivalent to show that $$1-\frac{1}{s^2-3s+3}\leq\sqrt{(s+t)(s+t-1)}.$$
It holds for $s\geq 3$ and $t\geq 1.$
%Since $t\geq (6-2\sqrt{5})s^2-(19-6\sqrt{5})s+17-6\sqrt{5}$, then $$s+t+1-\frac{s+t+1}{s^2-3s+3}\leq s+t-(5-2\sqrt{5})\leq\sqrt{(s+t)(s+t-1)}\quad(s+t\geq 5).$$

So we can assume that $H\subseteq G-\{v\}$ for any $v\in V(G)$. By Claim \ref{3.5} and Claim \ref{3.2}, we have $\omega(G_v)\leq s+t$ and $$x_v<\frac{1}{s+t+1}.$$
Since $\lambda(G)>\lambda(K_{s+t+1}^3)$ and $S_{2, s+t}$ is $\lambda$-perfect, then $S_{2, s+t}=\{v_1v_2u_1, v_1v_2u_2, ..., v_1v_2u_{s+t}\}\subseteq G$. By Claim \ref{3.4}, $G-\{v_1, v_2\}$ is $H$-free. Applying Claim \ref{3.6}, we have $\omega((G-\{v_2\})_{v_1})\leq s+t-1$ and $\omega((G-\{v_1\})_{v_2})\leq s+t-1$. Since $t\geq\frac{3}{2}s^2-\frac{11}{2}s+4$, then $$6\lambda(K_{s+t-1}^3)=\frac{(s+t-2)(s+t-3)}{(s+t-1)^2}\geq\frac{s^2-3s}{s^2-3s+2}>\pi_{\lambda}(H).$$ Assume the weight of $v_1$ and $v_2$ are $a_1$ and $a_2$ respectively, and $a_1+a_2=2a$, then
\begin{eqnarray*}
\lambda(G)&\leq &\frac{1}{6}\pi_{\lambda}(H)(1-2a)^3+a^2(1-2a)+2a\times\frac{1}{2}(1-\frac{1}{s+t-1})(1-2a)^2\\
&\leq &\frac{1}{6}\frac{s^2-3s}{(s-1)(s-2)}(1-2a)^3+a^2(1-2a)+2a\times\frac{1}{2}(1-\frac{1}{s+t-1})(1-2a)^2\\
&\leq &\lambda(K_{s+t-1}^3)(1-2a)^3+a^2(1-2a)+2a\times\frac{1}{2}(1-\frac{1}{s+t-1})(1-2a)^2\\
&=&\lambda(K_{s+t+1}^3, (a, a, \frac{1-2a}{s+t-1}, ..., \frac{1-2a}{s+t-1}))\\
&\leq&\lambda(K_{s+t+1}^3).
\end{eqnarray*}
\q

\section{Irrational Lagrangian densities}\label{irr}

In this section, we prove  Theorem \ref{P4}.

%\begin{defi}
%$S_i(n)$ is a 3-graph with vertex set [n] and edge set $\{v_1v_2v_3\in [n]^3| |\{v_1, v_2, v_3\}\cap [i]|>0\}$.
%\end{defi}

\begin{fact}\label{S2}
Let $V(S_2(n))$=[n] and $E(S_2(n))$=$\{12i| i\in [n]\setminus\{1, 2\}\}\cup\{ijk| i\in [2], j, k\in [n]\setminus\{1, 2\}\}$, then $\lambda(S_2(n))=\frac{\sqrt{3}}{18}$.
\end{fact}
\emph{Proof.} Let $\vec x=\{x_1, x_2, ..., x_n\}$ be an optimal weight vector on $S_2(n)$. By Fact \ref{symmetry}, we may assume that $x_1=x_2$, $x_3=x_4=...=x_n.$ Let $x_1=a$, then $x_3=\frac{1-2a}{n-2}$ and $a<\frac{1}{2}$. So
\begin{eqnarray*}
\lambda(S_2(n))&=&a^2(1-2a)+2a\bigg(\frac{1-2a}{n-2}\bigg)^2\binom{n-2}{2}\\
&\leq&2a^3-3a^2+a=f(a)\\
&\leq&f(\frac{3-\sqrt{3}}{6})=\frac{\sqrt{3}}{18}.
\end{eqnarray*}
\q\\
\emph{Proof of Theorem \ref{P4}.} Since $S_2(n)$ is $\mathscr{P}$-free, then by Fact \ref{S2} we have $\pi_{\lambda}(\mathscr{P})\geq3!\lambda(S_2(n))=\frac{\sqrt{3}}{3}$. For the upper bound, we assume that G is a $\mathscr{P}$-free dense 3-graph with vertex set $[n]$. It's sufficient to show that $\lambda(G)\leq\frac{\sqrt{3}}{18}$. Suppose on the contrary that $\lambda(G)>\frac{\sqrt{3}}{18}$, then $n>6$, since $\lambda(K_6^3)=\frac{5}{54}<\frac{\sqrt{3}}{18}$. By Proposition \ref{relationlt} and Theorem \ref{FL}, $\pi_{\lambda}(K_4^3)=\pi(K_4^3)<\frac{\sqrt{3}}{3}$. So $K_4^3\subseteq G$, without loss of generality, let \{1, 2, 3, 4\} form a $K_4^3$. Let $P=\{123, 124, 134, 234, 561, 562\}$.

\begin{lemma}\label{P1}
$P\nsubseteq G$.
\end{lemma}
\emph{Proof of Lemma \ref{P1}.} Suppose that $P\subseteq G$. We claim that for any $u_1, u_2\in V(G)\setminus\{1, 2, 3, 4, 5, 6\}$, $N(u_1, u_2)\subseteq\{1, 2\}$. If there exists $u_3\in V(G)\setminus\{1, 2, 3, 4, u_1, u_2\}$ such that $u_1u_2u_3\in E(G)$, then $\{1, 2, 3, 4, u_1, u_2, u_3\}$ forms a $P_1$. If there exist $u_3\in \{3, 4\}$ such that $u_1u_2u_3\in E(G)$, then $\{1, 2, 3, 4, 5, 6, u_1,\\ u_2, u_3\}$ forms a $P_2$. And $u_134\notin E(G)$ since otherwise $\{1, 2, 3, 4, 5, 6, u_1\}$ forms a $P_3$. And $u_135\notin E(G)$ since otherwise $\{1, 2, 3, 4, 5, 6, u_1\}$ forms a $P_4$, similarly $u_136, u_145, u_146\notin E(G)$. Since G is $P_1$-free, then every edge contains a vertex from $\{1, 2, 3, 4\}$. So we have shown that $G\subseteq S_2(n)$. By Fact \ref{S2}, $\lambda(G)\leq\frac{\sqrt{3}}{18}$, a contradiction.
\q

For any $v\in V(G)$ with weight $x_v$, we claim that $\omega(G_v)\geq 3$. Otherwise by Theorem \ref{MStheo} and Fact \ref{fact2}, we have $$3\lambda=\frac{\partial\lambda}{\partial x_v}\leq\frac{1}{2}(1-\frac{1}{2})(1-x_v)^2\leq\frac{1}{4}<\frac{\sqrt{3}}{6},$$ a contradiction. Let $x\in V(G)\setminus \{1, 2, 3, 4\}$. Since $\omega(G_x)\geq 3$, then $x$ is contained in a $K_4^{3-}$, denoted by $K_x$, and $d_{K_x}(x)=3$. Since G is $P_1$-free, then $|K_x\cap\{1, 2, 3, 4\}|$=2 or 3.

\emph{Case 1.} There exists $x\in V(G)\setminus\{1, 2, 3, 4\}$ such that $|K_x\cap\{1, 2, 3, 4\}|$=2.

Without loss of generality, let $K_x$=$\{1, 2, x, y\}$. Then $\{1, 2, 3, 4, x, y\}$ forms a $P$, by Lemma \ref{P1}, a contradiction.

\emph{Case 2.} For any $x\in V(G)\setminus\{1, 2, 3, 4\}$, $|K_x\cap\{1, 2, 3, 4\}|$=3.

Without loss of generality, let $K_x$=$\{1, 2, 3, x\}$ for some $x\in V(G)\setminus\{1, 2, 3, 4\}$. Then $K_x$ forms a $K_4^3$. We claim that for any $y \in V(G)\setminus\{1, 2, 3, 4, x\}$, $K_y$=$\{1, 2, 3, y\}$. If not, without loss of generality, let $K_y$=$\{2, 3, 4, y\}$. Then $K_x\cup\{y42, y43\}$ forms a $P$, by Lemma \ref{P1}, a contradiction. We claim that for any $v_1, v_2, v_3\in V(G)\setminus\{1, 2, 3\}$, $v_1v_2v_3\notin E(G)$ since otherwise there exist a $P_1$. And for any pair $v_1, v_2\in V(G)\setminus\{1, 2, 3\}$, we claim that $|N(v_1, v_2)\cap\{1, 2, 3\}|$=1. If $N(v_1, v_2)\cap\{1, 2, 3\}\geq2$, take $v_3\in V(G)\setminus\{1, 2, 3, v_1, v_2\}$, then $\{1, 2, 3, v_1, v_2, v_3\}$ forms a $P$, by Lemma \ref{P1}, a contradiction.

Assume the weight of $i\in V(G)$ is $x_i$ and $a$=$x_1+x_2+x_3$. Consider $\sum_{k=1}^3\frac{\partial\lambda}{\partial x_k}$. From the previous discussion we know that $x_ix_j$ only appears once in $\sum_{k=1}^3\frac{\partial\lambda}{\partial x_k}$ for $4\leq i, j\leq n$.

So
\begin{eqnarray*}
9\lambda&=&\sum_{k=1}^3\frac{\partial\lambda}{\partial x_k}\\
&\leq&\bigg(\sum_{i, j\in V(G)\setminus\{1, 2, 3\}}x_ix_j\bigg)+x_1x_2+x_1x_3+x_2x_3+2(x_1+x_2+x_3)\sum_{j\geq4}x_j\\
&=&\sum_{i, j}x_ix_j+(x_1+x_2+x_3)\sum_{j\geq4}x_j\\
&\leq&\frac{1}{2}+a(1-a)\leq\frac{3}{4}.
\end{eqnarray*}
Then $\lambda(G)\leq\frac{1}{12}<\frac{\sqrt{3}}{18}$, a contradiction.
\q

\end{document}